\documentclass[twoside,12pt]{article}

\usepackage{amsmath}
\usepackage{amssymb}
\usepackage{amscd}
\usepackage{amsthm}

\numberwithin{equation}{section}

\theoremstyle{plain}

\theoremstyle{remark}

\theoremstyle{definition}

\newcommand{\R}{\mathbb R}

\def\ra{\rightarrow}

\def\e{\emph}
\def\i{\infty}
\def\p{\partial}
\def\b{\begin}

\begin{document}

\title{
{Quasiisometries of   negatively curved   homogeneous manifolds associated with Heisenberg groups}}
\author{Xiangdong Xie\footnote{Partially supported by NSF grant DMS--1265735.}}
\date{  }

\maketitle

\begin{abstract}
We
  study  quasiisometries between  negatively curved homogeneous manifolds
  associated with diagonalizable derivations on Heisenberg   algebras.
  We classify these manifolds up to
quasiisometry,  and show that    all
 quasiisometries between  such manifolds
  (except
 when  they   are  complex
hyperbolic spaces)
  are almost similarities.
    We prove these results by studying the quasisymmetric
 maps on the ideal boundary of these manifolds.

\end{abstract}

{\bf{Keywords.}} quasiisometry, quasisymmetric map, negatively curved homogeneous manifolds,
    Heisenberg groups.



 {\small {\bf{Mathematics Subject
Classification (2010).}}  22E25,   30L10,  20F65.













\setcounter{section}{0}
  \setcounter{subsection}{0}

\section{Introduction}\label{s0}

In   this paper, we study quasiisometries between negatively curved homogeneous manifolds
  associated with Heisenberg groups.   We establish quasiisometric rigidity and quasiisometric classification results  for those manifolds
  associated with diagonalizable derivations.

Let $H_n$ be the $n$-th Heisenberg group  and $\mathcal{H}_n$ its Lie algebra.
    We shall identify  $H_n$ and $\mathcal{H}_n$   via the exponential map
  $\exp: {\mathcal H}_n\ra H_n$.
   Let  $A: \mathcal{H}_n  \ra  \mathcal{H}_n$  be a  derivation, that is,   $A$ is a
  linear  map satisfying  $A[X,Y]=[AX, Y]+[X, AY]$ for all
  $X, Y\in {\mathcal H}_n$.
    Define an action of $\R$ on $H_n$ by:
  $$t\cdot x =e^{tA} x\;\; \text{for}\;\;  x\in {\mathcal H}_n={ H}_n,  \;  t\in \R.$$
 Then one can form the semi-direct product
 $G_A=H_n\rtimes  \R$.    When the eigenvalues of $A$   have positive real parts, the group
 $G_A$ admits a left invariant Riemannian metric with negative sectional curvature  \cite{H}.
   In the case when $A$ is the standard derivation   with eigenvalues $1$ and $2$,
the manifold
 $G_A$  is isometric  to the complex hyperbolic space.

Assume $A: \mathcal{H}_n  \ra  \mathcal{H}_n$   is a diagonalizable   derivation. 
Suppose $A$ has positive eigenvalues
 $0<\alpha_1<\cdots <\alpha_k<\alpha_{k+1}$.  Let $U_i$ be the eigenspace associated with $\alpha_i$. Then we   have
 ${\mathcal H}_n=U_1\oplus \cdots \oplus U_k\oplus U_{k+1}$.  Every element $x\in {\mathcal H}_n$ can be written as
 $x=x_1+\cdots + x_k+x_{k+1}$ with $x_i\in U_i$.   By the above discussion,
    the group $G_A$ has a
left invariant Riemannian metric with negative sectional curvature.
  The ideal boundary $\p G_A$ can be naturally identified with  (the one point compactification of )   the      Heisenberg group $H_n$.  Fix a norm $|\cdot|$ on each    $U_i$.    
     The parabolic
   visual quasimetric  $d_A$ on $H_n={\mathcal H}_n$
 can be described as follows:
    $d_A(p, q)=||(-p)* q||_A$ for $p, q\in {\mathcal H}_n$,  where the norm $||\cdot||_A$ on ${\mathcal H}_n$ is   given by:
  $$||x_1+\cdots + x_k+x_{k+1}||_A=
\sum_{i=1}^{k+1} |x_i|^{\frac{1}{\alpha_i}}.$$

Similarly let $B: \mathcal{H}_n  \ra  \mathcal{H}_n$   be  a diagonalizable    derivation  with positive eigenvalues
 $0<\beta_1<\cdots< \beta_l<\beta_{l+1}$.   Let $W_i$ be the eigenspace of $\beta_i$.
The parabolic visual quasimetric  $d_B$ on 
 $H_n={\mathcal H}_n$     is similarly defined:
    $d_B(p, q)=||(-p)* q||_B$ for $p, q\in {\mathcal H}_n$,  where the norm $||\cdot||_B$ on ${\mathcal H}_n$ is  given by:
  $$||y_1+\cdots   + y_l+y_{l+1}||_B=
\sum_{i=1}^{l+1} |y_i|^{\frac{1}{\beta_i}}.$$

A map $f: X\ra Y$ between two quasimetric spaces is called an \e{almost
similarity}  if  there
 are
 constants $L\ge 1$ and $C\ge 0$ such that
   $L\cdot d(x_1,  x_2)-C\le d(f(x_1), f(x_2))\le L\cdot d(x_1,  x_2)+C$ for all
    $x_1,  x_2\in  X$   and $d(y, f(X))\le C$ for all $y\in Y$.

\b{Th}\label{main1}
  Let $A, B$   be diagonalizable    derivations  with positive eigenvalues,  and
  $G_A$ and $G_B$ the associated groups.   If   $k\ge 2$, then   every
  quasiisometry  $f: G_A  \ra G_B$ is an almost   similarity.

\end{Th}


\b{Th}\label{main2}
  Let $A, B$   be diagonalizable    derivations  with positive eigenvalues,  and
  $G_A$ and $G_B$ the associated groups.   Then
     $ G_A$  and    $  G_B$   are quasiisometric if and only if
   $k=l$, $\text{dim}(U_i)=\text{dim}(W_i)$,  and there is some $\lambda>0$ such that
 $\alpha_i=\lambda \beta_i$  for $1\le i\le k$.

\end{Th}

  Theorems \ref{main1}  and \ref{main2} generalize  the main results in \cite{SX}
     from the Euclidean group case to the    Heisenberg group case.
   The   general   case  for Euclidean groups were solved in   \cite{X}.
 The general case for the   Heisenberg groups   remains open.

 The   strategy   of  the  proof  is  the same as in \cite{SX}, that is, we study  quasisymmetric maps on the ideal boundary.
   In fact,   we shall prove  the following results for quasisymmetric maps.

\b{Th}\label{main3}
  Let $A, B$   be diagonalizable    derivations  with positive eigenvalues.
       If   $k\ge 2$  and    $F:   (H_n, d_A)  \ra (H_n, d_B) $  is a  quasisymmetry,  then
    $F$  is    biLipschitz  from $(H_n, d_A) $   to   $(H_n, d_B^{\frac{\beta_1}{\alpha_1}}) $.

\end{Th}

\b{Th}\label{main4}
  Let $A, B$   be diagonalizable    derivations  with positive eigenvalues.     Then
     $ (H_n, d_A)$  and $(H_n, d_B)$ are quasisymmetric  if and only if
   $k=l$, $\text{dim}(U_i)=\text{dim}(W_i)$,  and there is some $\lambda>0$ such that
 $\alpha_i=\lambda \beta_i$.

\end{Th}

The claims in   Theorem \ref{main1}   and Theorem \ref{main3} do not hold if $k=1$.   In this case,
 the manifold $G_A$ is the complex hyperbolic space and $(H_n, d_A)$ is   biLipschitz to the Heisenberg group with the Carnot metric.
 It is known that there are non-biLipschitz quasiconformal maps on the Heisenberg groups \cite{B}.
  Furthermore,  the claims in   Theorem \ref{main1}   and Theorem \ref{main3}   are equivalent, and so there are
 quasiisometries of the complex hyperbolic space that are not almost similarities.

Our results concern the quasiisometric rigidity and quasiisometric classification of negatively curved
solvable   Lie groups   $N\rtimes \R$.  The first result  in this area  is Pansu's  rigidity theorem \cite{P}   for the
  quarternionic hyperbolic spaces and Cayley plane.    The case $N=\R^n$ was solved in \cite{X}.
  In this paper we treat  the case $N=H_n$, but   only  for diagonalizable derivations.
   In \cite{X1}, \cite{X2} and \cite{X3},   we proved the   quasiisometric rigidity theorem for
$N\rtimes \R$  for many Carnot groups  $N$,  where  $\R$ acts on $N$ by the standard dilations on Carnot groups.
  All these  results  belong to  the larger project of
  quasiisometric rigidity and   quasiisometric  classification of focal hyperbolic groups \cite{C}.
  In this context,
Dymarz \cite{D}  recently     obtained  results similar to our Theorem \ref{main1} and
  Theorem \ref{main3}   for
  mixed type focal hyperbolic groups.

In   Section \ref{prelimi} we recall the definitions of various maps.
 In Section  \ref{derivation}
       we study the structure of diagonalizable   derivations on   $\mathcal H_n$.
In Section  \ref{metriconb}  we  define the homogeneous
manifolds associated with the Heisenberg groups,  and then study  the
visual quasimetric on
their ideal boundary.    
 In Section \ref{foliation}
   we show   that  every quasisymmetric map preserves a foliation.  In Section \ref{leaf} we show the restriction of a quasisymmetric map
 to a leaf is  biLipschitz.  Finally in Section \ref{proofs} we finish the proofs of the main theorems.

\noindent {\bf{Acknowledgment}}. {This work was
 initiated   while the author was attending the workshop
\lq\lq Interactions between analysis and geometry" at IPAM,  University of California at Los Angeles  from March to June 2013.  I would like to thank IPAM  for   financial support, excellent working conditions and
   conducive atmosphere.}

\section{Some basic definitions}\label{prelimi}

In this section we recall some basic definitions.

Let $K\ge 1$ and $C>0$. A bijection $F:X\ra Y$ between two
  quasimetric spaces is called a $(K,C)$-\e{quasisimilarity} if
\[
   \frac{C}{K}\, d(x,y)\le d(F(x), F(y))\le C\,K\,   d(x,y)
\]
for all $x,y \in X$.
   When $K=1$, we say $F$ is a \e{similarity}.
It is clear that a map is a quasisimilarity if and only if it is a
biLipschitz map. The point of using the notion of quasisimilarity is
that sometimes there is control on $K$ but not on $C$.

Let $L\ge 1$ and $A\ge 0$. A map $f: X\ra Y$ between two metric spaces is a
  $(L, A)$-quasiisometry if \newline
 (1)       for all  $ x_1, x_2\in X$:
     $$d(x_1, x_2)/L  -A\le d(f(x_1), f(x_2))\le L\cdot d(x_1, x_2)+A;$$
(2)     $d(y, f(X))\le A$ for all $y\in Y$. \newline
  The map $f$ is called a quasiisometry if it is a  $(L, A)$-quasiisometry for some $L\ge 1$, $A\ge 0$.

  Let $\eta: [0,\i)\ra [0,\i)$ be a homeomorphism.
    A bijection
$F:X\to Y$ between two   quasimetric spaces is
\e{$\eta$-quasisymmetric} if for all distinct triples $x,y,z\in X$,
we have
\[
   \frac{d(F(x), F(y))}{d(F(x), F(z))}\le \eta\left(\frac{d(x,y)}{d(x,z)}\right).
\]
If $F: X\ra Y$ is an $\eta$-quasisymmetry, then
  $F^{-1}:    Y\ra X$ is an $\eta_1$-quasisymmetry, where
$\eta_1(t)=(\eta^{-1}(t^{-1}))^{-1}$. See \cite{V}, Theorem 6.3.
    A map $F:  X\to Y$ is quasisymmetric if it is $\eta$-quasisymmetric
for some $\eta$.

Let $g: X_1\ra   X_2$ be a  bijection    between two
   quasimetric spaces  such that for any $p\in  X_1$,
     $d(x,p)\ra 0$ if and only if  $d(g(x), g(p))\ra 0$.
We define for every $x\in X_1$  and $r>0$,
\begin{align*}
   L_g(x,r)&=\sup\{d(g(x), g(x')):   d(x,x')\le r\},\\
   l_g(x,r)&=\inf\{d(g(x), g(x')):   d(x,x')\ge r\},
\end{align*}
and set
\[
   L_g(x)=\limsup_{r\ra 0}\frac{L_g(x,r)}{r}, \ \
   l_g(x)=\liminf_{r\ra 0}\frac{l_g(x,r)}{r}.
\]
  Then
\b{equation}\label{e5}
  L_{g^{-1}}(g(x))=\frac{1}{l_g(x)} \ \text{ and }\ l_{g^{-1}}(g(x))=\frac{1}{L_g(x)}
\end{equation}
for any $x\in X_1$. If $g$ is an $\eta$-quasisymmetry, then
$L_g(x,r)\le \eta(1)l_g(x, r)$ for all $x\in X_1$ and $r>0$. Hence
if in addition
\[
    \lim_{r\ra 0}\frac{L_g(x,r)}{r}\ \ {\text{or}} \ \ \lim_{r\ra 0}\frac{l_g(x,r)}{r}
\]
exists, then
\[
    0\le l_g(x)\le L_g(x)\le \eta(1) l_g(x)\le \infty.
\]

\section{Diagonalizable   derivations on   $\mathcal H_n$}\label{derivation}

In this section we   study  the structure of
  diagonalizable   derivations on   $\mathcal H_n$.

Let $A: \mathcal{H}_n\ra \mathcal{H}_n$  be a  diagonalizable   derivation.   Suppose  $A$    has   positive eigenvalues
  $0<\alpha_1<   \cdots < \alpha_k<\alpha_{k+1}$.  Let $U_i$ be the eigenspace of $\alpha_i$.
  Then we  have ${\mathcal H}_n=U_1\oplus \cdots \oplus U_k\oplus U_{k+1}$  and $A(v)=\alpha_i v$ for $v\in U_i$.
  For $x\in {\mathcal H}_n$, we write $x=x_1+\cdots + x_k+x_{k+1}$  with $x_i\in U_i$.

Fix a non-zero $e\in U_{k+1}$, and denote $m_i=\text{dim}(U_i)$.

\b{Le}\label{structure}
The following hold:\newline
(1) $m_{k+1}=1$ and $U_{k+1}=\mathcal Z(\mathcal H_n)$ is the center of ${\mathcal H}_n$;\newline
(2) $[U_i, U_j]=0$  if $i+j\not=k+1$;\newline
  (3)   $[U_i, U_{k+1-i}]\not=0$;\newline
(4)  $\alpha_i+\alpha_{k+1-i}=\alpha_{k+1}$  for all $1\le i\le k$;\newline
(5)  $m_i=m_{k+1-i}$  for $1\le i\le k$; \newline
  (6)    for   $i<(k+1)/2$,
   there exist  a  basis $e_1, \cdots, e_{m_i}$  for $U_i$  and   a  basis
 $\eta_1, \cdots, \eta_{m_i}$  for $U_{k+1-i}$  such that $[e_s, \eta_t]=\delta_{st} e$;\newline
   (7)  if $i=k+1-i$,  then   $m_i=2 k_i$ is even, and there is a basis $e_1, \eta_1, \cdots, e_{k_i}, \eta_{k_i}$ of
  $U_i$ such that  $[e_s, \eta_t]=\delta_{st} e$, $[e_s, e_t]=[\eta_s, \eta_t]=0$ for all $s,t$.

\end{Le}

\b{proof}
 (1).  Let  $X\in U_{k+1}$ be arbitrary.
Since $A$ is a derivation,
  for any   $ i$  and any  $Y\in U_i$,
we have  $A[X,Y]=[AX, Y]+[X, AY]=\alpha_{k+1}[X, Y]+\alpha_i[X, Y]=(\alpha_{k+1}+\alpha_i)[X, Y]$.
 As $\alpha_{k+1}$ is the largest eigenvalue of $A$ and $\alpha_i+\alpha_{k+1}>\alpha_{k+1}$,
  we must have $[X, Y]=0$.    This implies that $U_{k+1}\subset \mathcal Z(\mathcal H_n)$.
  Since $U_{k+1}$  is non-trivial and  $\mathcal Z(\mathcal H_n)$  has dimension $1$,     they  must   agree.

  (2), (3)   and (4).    We claim that for each $1\le  i\le k$, there exists a unique $j$ such that  $[U_i,  U_j]\not=0$.
  First of all, there exists at least one such $j$ since otherwise  $[U_i, {\mathcal H}_n]=0$  and so $U_i\subset U_{k+1}$.
   Now suppose there
  exist $j_1\not=j_2$  such that   $[U_i, U_{j_1}]\not=0$,   $[U_i, U_{j_2}]\not=0$.
   Then there are $X_1, X_2\in U_i$ and $Y_1\in U_{j_1}$,  $Y_2\in U_{j_2}$ such that
   $[X_1, Y_1]\not=0$  and $[X_2,  Y_2]\not=0$.
  Since $[{\mathcal H}_n,  {\mathcal H}_n]=U_{k+1}$,
   we have  $[X_1,  Y_1]=a_1 e$ and $[X_2,  Y_2]=a_2 e$ for some $a_1,  a_2\not=0$.
Since $A$ is a   derivation, we have
  $$A[X_1,Y_1]=[AX_1, Y_1]+[X_1,  AY_1]=\alpha_i [X_1, Y_1]+\alpha_{j_1} [X_1, Y_1]=(\alpha_i+\alpha_{j_1})[X_1, Y_1],$$
  which implies that $\alpha_i+\alpha_{j_1}=\alpha_{k+1}$.
  Similarly  by considering $A[X_2, Y_2]$   we
  obtain     $\alpha_i+\alpha_{j_2}=\alpha_{k+1}$.    It follows that  $\alpha_{j_1}=\alpha_{j_2}$, a contradiction.

We shall denote by $j_i$ the unique $j$  such that  $[U_i,  U_j]\not=0$.
  The preceding paragraph shows  that
  $\alpha_s+\alpha_{j_s}=\alpha_{k+1}$  for all $s$.
  Since $\alpha_1<\alpha_2<\cdots <\alpha_k$,    we see that $s<t$ implies  $j_s>j_t$.
  It is now   easy  to see that
  $j_i=k+1-i$. Hence (2), (3)   and (4)   hold.

 (5).  Suppose
 $m_i>m_{k+1-i}$  for some $i$.    Let $L(U_{k+1-i},  U_{k+1})$ be the vector space of linear maps from $U_{k+1-i}$  to $U_{k+1}$.
  Define   a  linear map $g:  U_i\ra   L(U_{k+1-i},  U_{k+1})$
 by   $g(X)=\text{ad}(X)|_{U_{k+1-i}}$,
  where $ad(X):   {\mathcal H}_n\ra  {\mathcal H}_n$ is given by
  $ad(X)Y=[X, Y]$  for $Y\in {\mathcal H}_n$.       Since $U_{k+1}$ is $1$-dimensional  and
 $\text{dim}(U_i)>\text{dim}(U_{k+1-i})$,     the kernel of $g$   is  non-trivial.
Hence there is some $X\in U_i\backslash\{0\}$  such that $[X, U_{k+1-i}]=0$.  Now it follows from (2) that $[X, {\mathcal H}_n]=0$,
  which is impossible in ${\mathcal H}_n$.  Therefore
 $m_i=m_{k+1-i}$  for all $i$.

  (6).    
  Let  $e_1\in U_i$ be a nonzero vector.    Then there is some $\eta_1\in U_{k+1-i}$  such that
 $[e_1, \eta_1]\not=0$. After   multiplying $\eta_1$  by a  nonzero constant, we may assume
$[e_1, \eta_1]=e$.
  Now let $e_2\in U_i\cap \text{ker}(\text{ad} (\eta_1))$ be a nonzero vector  and  as above pick
  $\eta_2\in U_{k+1-i}\cap \text{ker}(\text{ad} (e_1))$  such that
 $[e_2, \eta_2]=e$.   Inductively  we pick
 $$e_s\in U_i\cap \text{ker}(\text{ad} (\eta_1))\cap\cdots \cap \text{ker}(\text{ad} (\eta_{s-1}))$$
   and
  $$\eta_s\in  U_{k+1-i}\cap \text{ker}(\text{ad} (e_1))\cap\cdots \cap \text{ker}(\text{ad} (e_{s-1}))$$
    such that
  $[e_s, \eta_s]=e$.  By the choice we have  $[e_s, \eta_t]=\delta_{st} \, e$  for all $s, t$.

(7).  The proof is similar to that of (6).  First pick $e_1, \eta_1\in U_i$ such that $[e_1, \eta_1]=e$.
     Then  inductively pick
 $$e_s, \eta_s\in U_i\cap \bigcap_{t=1}^{s-1} \left(\text{ker}(\text{ad} (e_t))\cap \text{ker}(\text{ad} (\eta_{t}))\right)$$
  such that $[e_s, \eta_s]=e$.   In this way we get a basis satisfying all the conditions  in (7).

\end{proof}

\section{Quasimetric on the ideal boundary}\label{metriconb}

  The goal of  this Section   is to show
that the quasimetric $d_A$ defined in the Introduction is biLipschitz equivalent 
  with    a metric when the smallest  eigenvalue  of $A$  is at least $1$,
    see Lemma \ref{bilip2d}.   
     This result will be needed in Section \ref{foliation}  
for the application of Tyson's theorem (Theorem 1.4, \cite{T}): 
 Tyson's theorem does not apply to general quasimetric spaces.


Let $H_n$ be the $n$-th Heisenberg group and $\mathcal{H}_n$   its   Lie algebra.
  If  we identify  $\mathcal{H}_n$    with $\R^{2n}\times  \R=\R ^{2n+1}$,
  and if $e_i$, $1\le i\le 2n+1$ denote the standard basis of  $\R^{2n+1}$,
  then the only non-trivial Lie bracket relations are
 $[e_{2i-1}, e_{2i}]=e_{2n+1}$,  $1\le i\le n$.
  We shall identify $H_n$ and   $\mathcal{H}_n$  via the exponential map  and
 the group operation on $H_n$ shall be given by the BCH formula:
   $$X*Y=X+Y+\frac{1}{2}[X, Y]\;\;\; {\text for}\; \;X, Y\in {\mathcal H}_n.$$

  Let $A: \mathcal{H}_n\ra \mathcal{H}_n$  be a derivation, i.e.,  a linear map such that
 $$A[X, Y]=[AX, Y]+[X, AY]$$
 for all $X, Y\in \mathcal{H}_n$.
  Then one can define an action of $\R$  on $H_n$:
$$\R\times H_n \rightarrow   H_n$$
             $$(t, x)\rightarrow e^{tA} x.$$
             We denote the corresponding semi-direct product by
             $G_A=H_n\rtimes_A \R$.  Then $G_A$ is a solvable Lie
             group.
             Recall that the group operation in $G_A$ is given by:
             $$(g, t_1)\cdot (h, t_2)=(g* e^{t_1 A} h, \;   t_1
             + t_2).$$


By  Heintze's result (\cite{H}),  if
the eigenvalues of $A$ have positive real parts, then
 there is a left invariant Riemannian metric  on
  $G_A$
  with negative sectional curvature.
 Since any two left invariant Riemannian metrics are biLipschitz equivalent,
  $G_A$ is Gromov hyperbolic with any left invariant Riemannian metric.
   In this paper we only consider the  case
     when $A$ is  diagonalizable.

  Let
$A: \mathcal{H}_n  \ra  \mathcal{H}_n$   be  a diagonalizable    derivation  with positive eigenvalues.
Denote by  $0<\alpha_1<   \cdots < \alpha_k<\alpha_{k+1}$ the   eigenvalues  of
  $A$,  and   
    $U_i$  the eigenspace of $\alpha_i$.
  Then   $U_{k+1}=\mathcal Z(\mathcal H_n)$   and
${\mathcal H}_n=U_1\oplus \cdots \oplus U_k\oplus  U_{k+1}$.  
Fix a nonzero $e\in U_{k+1}$.  We choose  a  vector space basis $\mathcal B_i$ of $U_i$  for $1\le i\le n$  supplied by Lemma \ref{structure} (6), (7).
  Let  $\mathcal B=\cup_{i=1}^k \mathcal B_i\cup \{e\}$.
  On $T_e G_A$    we choose the inner product  such that $\mathcal B\cup \{\frac{\partial}{\partial t}\}$ is orthonormal.
  Let $g$ be the associated
left invariant Riemannian   metric  on $G_A$.



For each $g\in H_n$, the map $\gamma_g: \R\ra G_A$,
$\gamma_g(t)=(g,t)$  is a geodesic. We call such a geodesic a
vertical geodesic.  It can be checked that all vertical geodesics
are asymptotic as $t\ra +\infty$. Hence they define a point $\xi_0$
in the ideal boundary $\p G_A$.
    The sets $H_n\times\{t\}$  ($t\in \R$)
are horospheres centered at $\xi_0$, and $b:  G_A\ra \R$,   $b(x,t)=t$    is
 a
 Busemann   function
  associated  to   $\xi_0$.

Each geodesic ray in $G_A$ is  asymptotic to either  an upward
oriented vertical geodesic or a downward oriented vertical geodesic.
The upward oriented vertical geodesics are asymptotic to $\xi_0$ and
the downward oriented vertical  geodesics are in 1-to-1
correspondence with $H_n$. Hence $\p G_A\backslash\{\xi_0\}$ can be
naturally identified with $H_n$.

On  $T_e H_n=\mathcal H_n$ we fix the inner product such that $\mathcal B$ is orthonormal.
   Let $|\cdot|$ be the norm on $\mathcal H_n=(U_1\oplus \cdots \oplus U_k)\oplus U_{k+1}$
induced by this inner product.
For   $x=x_1+\cdots +x_k+x_{k+1}$ with $x_i\in U_i$,
 define
 $$||x||_A=\sum_{i=1}^{k+1}|x_i|^{\frac{1}{\alpha_i}}$$
 and
 $$||x||=|x_{k+1}|^{\frac{1}{2}}+\sum_{i=1}^{k}|x_i|.$$
For any $p,q \in \mathcal H_n$, let
$d_A(p,q)=||(-p)*q||_A$  and $d_0(p,q)=||(-p)*q||$.
  Notice that  $f(x)=a^x$ is a  non-increasing function  if  $0\le a\le 1$.
Suppose the smallest eigenvalue $\alpha_1$ of $A$ satisfies
 $\alpha_1\ge 1$.   Then
       $||x||_A\ge ||x||$
   whenever $||x||_A\le 1$.
Hence $d_A(p,q)\ge d_0(p,q)$ whenever $d_A(p,q)\le 1$.

Consider the grading  $\mathcal H_n=V_A\oplus \mathcal Z(\mathcal H_n)$,  where
  $V_A=U_1\oplus \cdots\oplus   U_k$.    For $t\in \R$, let
 $\delta_t:  \mathcal H_n \ra \mathcal H_n$  be the usual dilation on
 Heisenberg groups given by $\delta_t(v_1)=e^t v_1$ for $v_1\in V_A$ and $\delta_t(v_2)=e^{2t} v_2$ for $v_2\in \mathcal Z(\mathcal H_n)$.    It is well known that
  $d_0$ is biLipschitz equivalent with  every Carnot metric on $H_n$ associated with the
   above  grading.    We fix such a Carnot metric $d_C$ on $H_n$.  Then there is a constant $L\ge 1$ such that
  $$\frac{1}{L}\cdot d_0(p,q)\le d_C(p,q)\le L\cdot d_0(p,q)$$
 for any $p,q\in H_n$.

It follows from the definition of $||\cdot||_A$ and $d_A$ that 
\b{equation}\label{e1}
  d_A(g*x,   g*y)=d_A(x,y)\;\;\;  \text{for all}\;\;\; x,y, g\in
   H_n
\end{equation}
  and $||e^{tA}x||_A=e^t\cdot ||x||_A$ for  all  $x\in H_n$ and all $t\in \R$.  
 Since $e^{tA}$ is  an automorphism of  $H_n$, we have
  \b{equation}\label{e2}
  d_A(e^{tA}x, e^{tA}y)=e^t\cdot   d_A(x,y)\;\;\; \text{for all}\;\;\; x,y\in
       H_n  \;\;\;\text{and all} \;\;\;t\in \R.
\end{equation}



\b{Le}\label{bilip2d}
   If the smallest eigenvalue $\alpha_1$ of $A$   satisfies $\alpha_1\ge 1$,  then
  the quasimetric $d_A$ is biLipschitz equivalent   with  a  metric on $H_n$.
\end{Le}

\b{proof}
For any two points $p,q\in    H_n$, define
 $$\tilde d_A(p,q)=\inf\{\sum_{i=1}^m d_A(p_{i-1}, p_i): m\ge 1 \;\;\text{and}\;\;
 p_i\in H_n\;\; \text{with}\;\;  p_0=p, \;     p_m=q\}.$$
  We observe that $\tilde d_A$ also satisfies
 (\ref{e1}) and  (\ref{e2}).
   Let  $S=\{q\in H_n:  d_A(0, q)=1\}$  be the unit sphere with respect to $d_A$. \newline
{\bf{Claim:}}  if $\alpha_1\ge 1$, then  there is some $c$ satisfying $0<c<\frac{1}{2}$   such that for all $q\in S$:
 $$c \le \tilde d_A(0,q)\le 1.$$

  Since  both   $d_A$ and $\tilde d_A$ satisfy  (\ref{e1}) and  (\ref{e2}), the claim implies
  $c\cdot  d_A(p,q)  \le \tilde d_A(p,q)\le d_A(p,q)$  for all $p,q\in H_n$.
   It follows  that $\tilde d_A$ is a metric on $H_n$ and that $d_A$ is biLipschitz equivalent   with    the metric $\tilde d_A$.
We next prove the claim.

 Clearly we have $\tilde d_A(p,q)\le d_A(p,q)$ for all $p,q\in H_n$.
So we only need to prove the first inequality.
 If $\tilde d_A(0,q)\ge \frac{1}{2}$, then we are done.
   Now assume $\tilde d_A(0,q)<\frac{1}{2}$.
Let   $p_0,  p_1, \cdots, p_m$  be a finite sequence of points in $H_n$ such that $p_0=0$, $p_m=q$ and
 $\sum_{i=1}^m d_A(p_{i-1}, p_i)<\frac{1}{2}$.
  Then  $d_A(p_{i-1}, p_i)<\frac{1}{2}$  for each $i$.
  Since $\alpha_1\ge 1$, we
have
  $d_A(p_{i-1}, p_i)\ge  d_0(p_{i-1}, p_i)$. It follows that
$$\sum_{i=1}^m d_A(p_{i-1}, p_i)\ge  \sum_{i=1}^m d_0(p_{i-1}, p_i)\ge \frac{1}{L}
 \sum_{i=1}^m d_C(p_{i-1}, p_i)\ge  \frac{1}{L}   d_C(0, q)\ge \frac{a}{L},   $$
  where  $a=\min\{d_C(0, q):  q\in S\}$.   Since  $d_C(0, \cdot)$ is  continuous  and  $S$ is a compact subset in $H_n$ disjoint from $0$,  we must have $a>0$.
 Now it is clear that the claim holds for $c=\min\{1/2, a/L\}$.


\end{proof}


\section{Quasisymmetric maps preserve a   foliation}\label{foliation}

In this Section we show that  quasisymmetric maps  $F: (H_n, d_A)\ra (H_n, d_B)$
  send a  foliation on $ (H_n, d_A)$ to a  foliation on $ (H_n, d_B)$.

 Let
$A: \mathcal{H}_n  \ra  \mathcal{H}_n$   be  a diagonalizable    derivation  with positive eigenvalues
  $0<\alpha_1<   \cdots < \alpha_k<\alpha_{k+1}$.
We will use the notation from  the previous section. 
   In particular, $\mathcal B$ is the basis of $\mathcal H_n$  constructed in the last section. 
Let $m$ be the Lesbegue measure on ${\mathcal H}_n$ with respect to this basis.  Then $m$ is invariant under left translations, as the  Jacobian matrix
  of the left translations with respect to  the basis have determinant $1$.    Furthermore, the automorphism
 $e^{At}$    has matrix representation given by a block diagonal matrix $[e^{\alpha_1 t}I_{m_1}, \cdots,  e^{\alpha_k t}I_{m_k},  e^{(\alpha_1+\alpha_k)t}I_1]$,
  where $I_m$ is the   $m\times m$ identity matrix.
  Lemma \ref{structure} (4), (5)   imply that  the determinant  of $e^{At}$ equals  $e^{t(n+1)(\alpha_1+\alpha_k)}$.
   Hence, for any metric ball  $B(x, r)$ in $(H_n, d_A)$  with radius $r=e^t$, we  have
  $m(B(x, r))=m(B(o, e^t))=m(e^{At}B(o, 1))=r^{(n+1)(\alpha_1+\alpha_k)}\cdot m(B(o,1))$.  In particular,
 $m$ is Ahlfors $Q$-regular with $Q=(n+1)(\alpha_1+\alpha_k)$.

 Set $V_A=U_1\oplus \cdots \oplus U_k$.
    Now   define a   quasimetric $D_A$  on $V_A$  by:
   $$D_A(x_1+\cdots +x_k, \,  y_1+\cdots   +y_k)=\sum_{i=1}^k |x_i-y_i|^{\frac{1}{\alpha_i}}.$$
   Let $\pi:  (\mathcal{H}_n,  d_A)  \ra (V_A,   D_A) $ be the natural projection given by $\pi(x+ z)=x$  for $x\in V_A$, $z\in \mathcal Z(\mathcal H_n)$.
  We observe that $\pi$  is a  $1$-Lipschitz  map.

When  $k\ge 2$, Lemma \ref{structure} (2) implies  $[U_1, U_1]=0$.
   So   $U_1$ is a Lie  subalgebra of ${\mathcal H}_n$.  We will   abuse notation and also use $U_1$ to denote the connected Lie subgroup
of  $H_n$ with Lie algebra $U_1$.

\b{Le}\label{rec}
   Suppose $k\ge 2$ and   $\alpha_1=1$.
  Then every rectifiable curve in $(H_n, d_A)$ is contained in some left coset of $U_1$.

\end{Le}

\b{proof}
Let  $\gamma: [0,1]\ra (H_n, d_A)$  be a rectifiable curve.  Since $\pi: (H_n,  d_A) \ra (V_A,  D_A)$
   is   $1$-Lipschitz, the curve $\pi\circ \gamma$   is  a  rectifiable  curve in $(V_A,    D_A)$.
   Since $\alpha_i>1$  for each $i\ge 2$,   there is no non-trivial rectifiable curve in
  $(U_i,\,  |\cdot|^{\frac{1}{\alpha_i}})$ for $i\ge 2$.     Hence there are
$x_i\in U_i$ for each $i\ge 2$ such that $\pi\circ \gamma$  lies in $U_1\times \{x_2\}\times \cdots \times \{x_k\}$  and so $\gamma$ lies in   $U_1\times \{x_2\}\times \cdots \times \{x_k\}\oplus \mathcal Z({\mathcal H}_n)$.
Define a metric $D'$ on $U_1\times \mathcal Z({\mathcal H}_n)$ as follows:
    $$D'((x_1, z),   (x'_1, z'))=|x'_1-x_1|+|z'-z|^{\frac{1}{1+\alpha_k}}.$$
  Let
  $$f: (U_1\times \{x_2\}\times \cdots \times \{x_k\}\oplus \mathcal Z({\mathcal H}_n),  d_A) \ra (U_1\times  \mathcal Z({\mathcal H}_n),  D') $$
   be
  defined   by
$f(x_1+ x_2+ \cdots+ x_k +x_{k+1})=(x_1,   x_{k+1}+\frac{1}{2}[x_1, x_k])$.  Then it is easy to check that
 $f$ is an isometry.
     In  $ (U_1\times  \mathcal Z({\mathcal H}_n),  D') $  the  rectifiable curves   lie in subsets of the form $U_1\times \{p\}$ with $p\in \mathcal Z({\mathcal H}_n)$.  It follows that the only rectifiable curves in $(H_n,  d_A)$ lie  in   subsets of the   form
 $$\left\{(x_1+   x_2+ \cdots+  x_k+    (p-\frac{1}{2}[x_1, x_k]))|  x_1\in U_1\right\}=(x_2+\cdots +x_k+  p)* U_1,  $$
   where $x_i\in U_i$, $2\le i\le k$  and $p\in \mathcal Z({\mathcal H}_n)$ are fixed.
  These subsets are exactly the  left cosets of $U_1$.

\end{proof}

   Now let
$B: \mathcal{H}_n\ra \mathcal{H}_n$  be another   diagonalizable   derivation  with   positive eigenvalues
  $0<\beta_1<   \beta_2<\cdots < \beta_l<\beta_{l+1}$.  Let $W_j$ be the eigenspace of $\beta_j$.
  Then we  have   $W_{l+1}=\mathcal Z({\mathcal H}_n)$  and
${\mathcal H}_n=V_B\oplus  W_{l+1}$,
  where   $V_B= W_1\oplus \cdots \oplus W_l$.  
       As in the case of   $||\cdot||_A$ and $d_A$,  we    fix a basis for $\mathcal H_n$ supplied by Lemma \ref{structure}  and define norm $||\cdot||_B$   and quasimetric $d_B$ on $\mathcal H_n$.  
  For $y\in {\mathcal H}_n$, we write $y=y_1+\cdots + y_l+y_{l+1}$  with $y_j\in   W_j$.
     Then  
 $$||y_1+\cdots +y_l+y_{l+1}||_B=
\sum_{j=1}^{l+1} |y_j|^{\frac{1}{\beta_j}}.$$
 The visual  quasimetric   $d_B$   on $H_n={\mathcal H}_n$ is given by:
  $d_B(p, q)=||(-p)*q||_B$.



\b{Prop}\label{prefo}
Let
$F:  (H_n, d_A)\ra (H_n, d_B)$  be an $\eta$-quasisymmetric  map for some $\eta$.    Suppose $k\ge 2$.
  Then   $l\ge 2$,
$\text{dim}(U_1)=\text{dim}(W_1)$  and $F$ maps  left cosets of $U_1$  to left cosets of $W_1$.

\end{Prop}

\b{proof}
 By replacing $d_A$ and $d_B$ with suitable powers, we may assume $\alpha_1+\alpha_k=\beta_1+\beta_l$
  and $\min\{\alpha_1, \beta_1\}=1$.   Then
 $(H_n, d_A)$ and $(H_n, d_B)$ have the same Hausdorff dimension $Q=(n+1)(\alpha_1+\alpha_k)$.
  By  considering $F^{-1}$ instead of  $F$ if necessary  we may assume $\alpha_1=1$.
 So $\beta_1\ge 1$.  By Lemma \ref{bilip2d},
$(H_n, d_A)$ and $(H_n, d_B)$   are biLipschitz equivalent to metric  spaces.

We claim that $\beta_1=1$.  Suppose $\beta_1>1$.    Then  Lemma   \ref{rec}
and its proof show that there is no non-trivial rectifiable curve
in $(H_n, d_B)$.  In particular, every curve family in
$(H_n, d_B)$  has $Q$-modulus $0$.
 On the other hand, fix a nonzero vector $v\in U_1$.     Since $\alpha_1=1$, the definition of $d_A$ implies that
    the left translates of the segment $\sigma:=\{tv:  t\in [0,1]\}$
   are rectifiable. Let $U\subset {\mathcal H}_n$ be a hyperplane  transversal to  the direction $v$.
   By a classical calculation, the   family of curves $\Gamma:=\{g\cdot \sigma:  g\in U\}$
  has positive $Q$-modulus.    Since
  $(H_n, d_A)$ and $(H_n, d_B)$   are quasisymmetric and have the same Hausdorff dimension $Q>1$,
   by    Tyson's theorem (Theorem 1.4 in \cite{T}),   $F(\Gamma)$ also has positive   $Q$-modulus, contradicting the above  observation. Hence $\beta_1=1$.
  Then we also have $\alpha_k=\beta_l$.  We remark that Tyson's theorem holds only for metric spaces.  Since  the quasimetric spaces 
 $(H_n, d_A)$ and $(H_n, d_B)$     are biLipschitz equivalent with  metric  spaces,  we can still apply Tyson's theorem.

We next claim that for any left translate $g\cdot \sigma$ of the segment $\sigma$ as above,
 the image $F(g\cdot \sigma)$ lies in a left coset of $W_1$.  Since   any
     two points in a left  coset of $U_1$ can be joined by a segment of the form
 $g\cdot \sigma$,   the claim implies  that $F$ maps every left coset of $U_1$ into a left coset of $W_1$.   The same argument applied to $F^{-1}$ shows $F^{-1}$ maps left cosets of $W_1$ into left cosets of $U_1$. Hence  the image of a left coset of $U_1$ under $F$ is a left coset of $W_1$.  Next we prove the claim.

  Suppose   $F(g* \sigma)$ is not contained in any left coset of $W_1$.   By continuity of $F$, there is an open subset $U$ containing $g$ such that
 for any $g'\in U$, the image $F(g'* \sigma)$ also does not lie in any left coset of $W_1$.  By  Lemma   \ref{rec},
 $F(g'* \sigma)$  is not rectifiable.   So  the $Q$-modulus of the curve family $F(\Gamma)$ is $0$, where $\Gamma=\{g'* \gamma: g'\in U\}$.
    On the   other hand, as indicated above, the $Q$-modulus of $\Gamma$ is positive, contradicting   Tyson's theorem.  Hence the  claim holds.

\end{proof}

\section{Restriction to a leaf}\label{leaf}

In this Section we show that the restriction of a quasisymmetric map
  $F: (H_n, d_A)\ra (H_n, d_B)$
to a left coset of $U_1$ is
   a  quasisimilarity.

For the rest of this Section, let
$A, B: \mathcal{H}_n  \ra  \mathcal{H}_n$   be   diagonalizable    derivations  with positive eigenvalues.
Denote by  $0<\alpha_1<   \cdots < \alpha_k<\alpha_{k+1}$ the   eigenvalues  of
  $A$,  and   
    $U_i$  the eigenspace of $\alpha_i$.
  Then we  have   $U_{k+1}=\mathcal Z({\mathcal H}_n)$ and
${\mathcal H}_n=U_1\oplus \cdots \oplus U_k\oplus \mathcal Z({\mathcal H}_n)$.  
  Similarly   let
  $0<\beta_1<   \beta_2<\cdots < \beta_l<\beta_{l+1}$  be the  eigenvalues
of   $B$,
    $W_j$ be the eigenspace of $\beta_j$.
  Then we  have  $W_{l+1}=\mathcal Z({\mathcal H}_n)$ and
 ${\mathcal H}_n=W_1\oplus \cdots \oplus W_l\oplus \mathcal Z({\mathcal H}_n)$.
Without loss of generality, we may assume $\alpha_1=\beta_1=1$.

Fix a non-zero $e\in \mathcal Z({\mathcal H}_n)$.   We choose norms on $\mathcal H_n$ and define  quasimetrics 
   $d_A$ and $d_B$    on $H_n$  as    in 
  Section \ref{metriconb}.  In particular, we have $|e|=1$.  
  Note that  the Cauchy-Schwartz inequality implies $|[x_1, x_k]|\le |x_1|\cdot |x_k|$ for any
 $x_1\in U_1$,  $x_k\in U_k$.


  Now let
 $F: (H_n, d_A)\ra (H_n, d_B)$ be   an $\eta$-quasisymmetric map for some $\eta$.
By Proposition \ref{prefo},  $F$ maps  left cosets of $U_1$ to left cosets of $W_1$.
  Recall that $F^{-1}$ is $\eta_1$-quasisymmetric, where $\eta_1(t)=(\eta^{-1}(t^{-1}))^{-1}$.
 Without loss of generality we may assume $\eta(1)\ge 1$. Then we also have $\eta_1(1)\ge 1$.


The proof of the following Lemma is similar to that of   Lemma 5.1 in \cite{X1}.
  But the calculations are different.

\b{Le}\label{key}
Let $L$  be   a left coset of $U_1$ and  denote $L'=F(L)$.
   Suppose $p,q\in L$ are such that $l_F(p)>  C_1\cdot  L_F(q)$  with $C_1=102\cdot \eta_1(1)$.
  Let $c:\R\ra L$ be the parametrization of the line through $p$, $q$   such that 
 $c(0)=p$ and $c(1)=q$.  Then
$$ \text{L}_F(c(\lambda))\le  2{(\eta_1(1))^2}\left(\frac{2}{|\lambda|}\right)^{\frac{1}{1+\alpha_k}}\cdot  \text{L}_F(q)$$
 for all $|\lambda|\ge 1$;  in particular,
$L_F(c(\lambda))\le   3{(\eta_1(1))^2} L_F(q)$.


\end{Le}

\b{proof}    Denote    $p'=F(p)$  and $q'=F(q)$.
  The assumption implies   $l_{F^{-1}}(q')>  C_1\cdot  L_{F^{-1}}(p')$.
Let $\{r_j\}$ be an arbitrary sequence of positive reals such that $r_j\ra 0$.
  Then
$$\liminf_{j\ra \infty}\frac{l_{F^{-1}}(q', r_j)}{r_j}> C_1 \cdot \limsup_{j\ra \infty}\frac{L_{F^{-1}}(p', r_j)}{r_j}.$$
    We shall look at the image of  the left coset $r^{1+\beta_l}_j e+L'$ of  $W_1$  under $F^{-1}$.
  Recall  $e\in \mathcal Z({\mathcal H}_n)$ is a fixed element with $|e|=1$.
 By  Lemma \ref{prefo},  $L_j:=F^{-1}(r^{1+\beta_l}_j e+L')$ is a left coset of $U_1$.

Denote $p'_j=r^{1+\beta_l}_j e+p'$ and $q'_j=r^{1+\beta_l}_j e+q'$.  Notice that $d_B(p', p'_j)=r_j$ and $d_B(q', q'_j)=r_j$.
    So   we have
  $$\frac{d_A(q,  F^{-1}(q'_j))}{r_j}\ge \frac{l_{F^{-1}}(q',   r_j)}{r_j}$$
  and
   $$\frac{d_A(p,  F^{-1}(p'_j))}{r_j}\le \frac{L_{F^{-1}}(p',   r_j)}{r_j}.$$
  Let $p_j,  q_j\in   L_j$ be a point on $L_j$ nearest to $p$ and $q$ respectively.
  Notice that $p'_j$ is  the point on  $r^{1+\beta_l}_j e+L'$  nearest to  $p'$. 
     Hence  $d_B(p',  p'_j)\le d_B(p',  F(p_j))$. Since 
$F^{-1}$ is $\eta_1$-quasisymmetric, we have
   $d_A(p,  F^{-1}(p'_j))\le \eta_1(1) d_A(p,  p_j)$.  Similarly we have
 $d_A(q,  F^{-1}(q'_j))\le \eta_1(1) d_A(q,  q_j)$.
 It follows that
\b{equation}\label{e9}
\frac{d_A(q,  q_j)}{r_j}\ge \frac{1}{\eta_1(1)}\cdot \frac{d_A(q,  F^{-1}(q'_j))}{r_j}\ge
\frac{1}{\eta_1(1)}\cdot    \frac{l_{F^{-1}}(q',   r_j)}{r_j}
\end{equation}
  and
$$\frac{d_A(p,  p_j)}{r_j}\le
\frac{d_A(p,  F^{-1}(p'_j))}{r_j}\le \frac{L_{F^{-1}}(p',   r_j)}{r_j}.$$
      Therefore
  \b{align*}
\liminf_{j\ra \infty} \frac{d_A(q,  q_j)}{r_j}\ge \liminf_{j\ra \infty}   \frac{1}{\eta_1(1)}\cdot  \frac{l_{F^{-1}}(q',   r_j)}{r_j}
  &\ge    \frac{C_1}{\eta_1(1)}\cdot  \limsup_{j\ra \infty} \frac{L_{F^{-1}}(p',   r_j)}{r_j}\\
  & \ge    {102}\cdot  \limsup_{j\ra \infty}
    \frac{d_A(p,  p_j)}{r_j}.
\end{align*}
    Hence
$$\frac{d_A(p, p_j)}{d_A(q, q_j)}\le \frac{1}{101}$$
  for all sufficiently large $j$.

Next we shall look at $d_A(p, p_j)$  and $d_A(q, q_j)$.

Notice that
 $L=q*U_1=q*\{t: t\in U_1\}$.    Write $q_j=q*(\tilde{x}_1+\tilde{x}_2+\cdots  +\tilde{x}_k+ \tilde{z})$
  with $\tilde{x}_i\in U_i$ and $\tilde{z}\in \mathcal Z({\mathcal H}_n)$.   Although the $\tilde{x}_i$'s  and  $\tilde{z}$ depend on
 $r_j$,  we shall suppress the dependence to simplify the notation.
  Then $L_j=q_j*U_1=q_j*\{t:  t\in U_1\}$.
  An arbitrary point on $L_j$ has the form
  $$q_j*t'=q*(\tilde{x}_1+\tilde{x}_2+\cdots  +\tilde{x}_k+ \tilde{z})*t'
=q*\left((t'+\tilde{x}_1)+\tilde{x}_2+\cdots  +\tilde{x}_k+(\tilde{z}+\frac{1}{2}[\tilde{x}_k, t'])\right).$$
Since $q_j$ is   a   point on  $L_j$  nearest to $q$,
  we see that
 \b{align*}
  d_A(q, \, q_j*t')      &=||(-q)*q_j*t'||_A  \\
  &  =\left|\left|(t'+\tilde{x}_1)+\tilde{x}_2+\cdots  +\tilde{x}_k+   (\tilde{z}+\frac{1}{2}[\tilde{x}_k, t'])\right|\right|_A\\
  &  =|t'+\tilde{x}_1|+ |\tilde{x}_2|^{\frac{1}{\alpha_2}}+\cdots  +  |\tilde{x}_k|^{\frac{1}{\alpha_k}}+
\left|\tilde{z}+\frac{1}{2}[\tilde{x}_k, t']\right|^{\frac{1}{1+\alpha_k}}
\end{align*}
    achieves   minimal when $t'=0$.

Now write $p=q*t_0$  and $p_j=q_j*t_j\in    L_j$   for some $t_0,  t_j\in U_1$.   Then we have
\b{align*}
(-p)*p_j & =(-t_0)*(\tilde{x}_1+\tilde{x}_2+\cdots  +\tilde{x}_k+   \tilde{z})*t_j\\
 &  =(\tilde{x}_1+t_j-t_0)+\tilde{x}_2+\cdots  +\tilde{x}_k+  (\tilde{z}+\frac{1}{2}[\tilde{x}_k,  t_0+t_j]).
  \end{align*}
    Hence
  $$d_A(p, p_j)=|\tilde{x}_1+t_j-t_0|+|\tilde{x}_2|^{\frac{1}{\alpha_2}}+\cdots  +  |{\tilde x}_k|^{\frac{1}{\alpha_k}}+
   \left|\tilde{z}+\frac{1}{2}[\tilde{x}_k, t_0+t_j]\right|^{\frac{1}{1+\alpha_k}}.$$
   To simplify notation set $a=|\tilde{x}_2|^{\frac{1}{\alpha_2}}+\cdots  +  |\tilde{x}_k|^{\frac{1}{\alpha_k}}$.
We have
  $$d_A(q, q_j)=|\tilde{x}_1|+a+|\tilde{z}|^{\frac{1}{1+\alpha_k}}$$
  and
 $$d_A(p, p_j)=|\tilde{x}_1+t_j-t_0|+a+\left|\tilde{z}+\frac{1}{2}[\tilde{x}_k, t_0+t_j]\right|^{\frac{1}{1+\alpha_k}}.$$
 Since  $\frac{d_A(p, p_j)}{d_A(q, q_j)}\le \frac{1}{101}$,
   we  have $100a\le |\tilde{x}_1|+|\tilde{z}|^{\frac{1}{1+\alpha_k}}$.

\noindent
{\bf{Claim.}}
     $5|\tilde{x}_1|\le |\tilde{z}|^{\frac{1}{1+\alpha_k}}$.

  Suppose the contrary.  
  Then $|\tilde{z}|^{\frac{1}{1+\alpha_k}}< 5 |\tilde{x}_1|$.
Now
$$|\tilde{x}_k|^{\frac{1}{\alpha_k}}\le a\le \frac{|\tilde{x}_1|+|\tilde{z}|^{\frac{1}{1+\alpha_k}}}{100}\le
\frac{(5+1)|\tilde{x}_1|}{100}\le \frac{1}{10}|\tilde{x}_1|,
$$
  so
  $$|\tilde{x}_k|\le  \frac{1}{10^{\alpha_k}}|\tilde{x}_1|^{\alpha_k}.$$
It follows that $|[\tilde{x}_k, \tilde{x}_1]|\le  \frac{1}{10^{\alpha_k}}|\tilde{x}_1|^{1+\alpha_k}.$    Now
 \b{align*}
 d_A(q,\, q_j*(-\tilde{x}_1))&  =a+\left|\tilde{z}-\frac{1}{2}[\tilde{x}_k, \tilde{x}_1]\right|^{\frac{1}{1+\alpha_k}}\\
  &  \le
  a+|\tilde{z}|^{\frac{1}{1+\alpha_k}}  +  \left(\frac{1}{2}|[\tilde{x}_k, \tilde{x}_1]|\right)^{\frac{1}{1+\alpha_k}}\\
  &  \le a+|\tilde{z}|^{\frac{1}{1+\alpha_k}}  +\left(\frac{1}{2\cdot  10^{\alpha_k}}\right)^{\frac{1}{1+\alpha_k}} |\tilde{x}_1|\\
& < d_A(q, q_j),
\end{align*}
  contradicting the fact that $q_j$ is a point on $L_j$  nearest to $q$.  Hence the claim holds.

The above claim together with the estimate on $a$  implies
\b{equation}\label{e10}
d_A(q, q_j)\le2 |\tilde{z}|^{\frac{1}{1+\alpha_k}}.
\end{equation}
 So
 \b{equation}\label{e11}
 |\tilde{x}_k|^{\frac{1}{\alpha_k}}\le a \le d_A(p, p_j)\le \frac{1}{101 } d_A(q, q_j)
  \le  \frac{2}{101 }  |\tilde{z}|^{\frac{1}{1+\alpha_k}}.
\end{equation}

Now let $u=\tilde{x}_1+t_j-t_0$.
 Then $|u|\le d_A(p, p_j)\le   \frac{2}{101 }
     |\tilde{z}|^{\frac{1}{1+\alpha_k}}$.
  It follows that  $|[\tilde{x}_k,  u]|\le \left(\frac{2}{101}\right)^{1+\alpha_k}  |\tilde{z}| $. Similarly
  \b{equation}\label{e12}
|[\tilde{x}_k,  \tilde{x}_1]|\le    \frac{1}{5}\cdot \left(\frac{2}{101}\right)^{\alpha_k}  |\tilde{z}|.
\end{equation}
  On the other hand,
  $$\left|\tilde{z}+[\tilde{x}_k, t_0]+\frac{1}{2}[\tilde{x}_k, u]-  \frac{1}{2}[\tilde{x}_k, \tilde{x}_1]\right|^{\frac{1}{1+\alpha_k}}
=\left|\tilde{z}+\frac{1}{2}[\tilde{x}_k, t_0+t_j]\right|^{\frac{1}{1+\alpha_k}}\le d_A(p, p_j)\le \frac{2}{101 }  |\tilde{z}|^{\frac{1}{1+\alpha_k}} .$$
  Now the  triangle inequality implies   
 \b{equation}\label{e13}
|\tilde{z}+[\tilde{x}_k, t_0]|\le  \frac{1}{25} |\tilde{z}|\;\;\text{and}\;\;
   \frac{24}{25} |\tilde{z}|    \le  |[\tilde{x}_k, t_0]|\le  \frac{26}{25} |\tilde{z}|.
\end{equation}

For $\lambda\in \R$, denote   $w_\lambda=q*(\lambda t_0)\in L$.
Let $w_{\lambda, j}\in   L_j$ be a point on $L_j$ nearest to
   $w_\lambda$.  Then
$w_{\lambda, j}=q*(\tilde{x}_1+\tilde{x}_2+\cdots  +\tilde{x}_k+  \tilde{z})*t_\lambda$   for some $t_\lambda\in U_1$.
  By a calculation
   similar to that   of $d_A(p, p_j)$,  we obtain
$$d_A(w_\lambda, w_{\lambda,j})=|\tilde{x}_1+t_\lambda-\lambda t_0|+a+
\left|\tilde{z}+\frac{1}{2}[\tilde{x}_k, \lambda t_0+t_\lambda]\right|^{\frac{1}{1+\alpha_k}}.$$
Let $w=\tilde{x}_1+t_\lambda-\lambda t_0$.  Then $\lambda t_0+t_\lambda=2 \lambda t_0+w-\tilde{x}_1$.

  Now  suppose $|\lambda-1|\ge 1$.
 If  $|w|\ge \sqrt{|\lambda-1|}\cdot |\tilde{z}|^{\frac{1}{1+\alpha_k}}$,  then   by      (\ref{e10})
   $$d_A(w_\lambda, w_{\lambda,j})\ge  |w|\ge \frac{\sqrt{|\lambda-1|}}{2} d_A(q, q_j).$$
 Assume now that
$|w|\le \sqrt{|\lambda-1|}\cdot  |\tilde{z}|^{\frac{1}{1+\alpha_k}}$.
  Then  by (\ref{e11}),
$|[\tilde{x}_k, w]|\le (\frac{2}{101})^{\alpha_k}     \sqrt{|\lambda-1|} \cdot   |\tilde z|$.    It   now follows   from (\ref{e13}),
  (\ref{e12}),   (\ref{e10})
  and the triangle inequality that
\b{align*}
d_A(w_\lambda, w_{\lambda,j})   &\ge \left|\tilde{z}+\frac{1}{2}[\tilde{x}_k, \lambda t_0+t_\lambda]\right|^{\frac{1}{1+\alpha_k}}\\
  &  = \left|\tilde{z}+\lambda [\tilde{x}_k,   t_0]+\frac{1}{2}[\tilde{x}_k, w]-\frac{1}{2}[\tilde{x}_k, \tilde{x}_1]
    \right|^{\frac{1}{1+\alpha_k}}\\
&  = \left|(\lambda-1) [\tilde{x}_k,   t_0]+
(\tilde{z}+ [\tilde{x}_k,   t_0])+\frac{1}{2}[\tilde{x}_k, w]-\frac{1}{2}[\tilde{x}_k, \tilde{x}_1]
    \right|^{\frac{1}{1+\alpha_k}}\\
&\ge \left(|\lambda-1|\cdot \frac{24}{25} |\tilde z| -\frac{1}{25}|\tilde z|-\frac{1}{2}\left(\frac{2}{101}\right)^{\alpha_k}\sqrt{|\lambda-1|}\cdot |\tilde z|
-\frac{1}{10}  \left(\frac{2}{101}\right)^{\alpha_k}  \cdot |\tilde z|    \right)^{\frac{1}{1+\alpha_k}}\\
    & \ge \left(\frac{|\lambda-1|}{2}\right)^{\frac{1}{1+\alpha_k}}\cdot |\tilde{z}|^{\frac{1}{1+\alpha_k}}\\
  &  \ge    \frac{1}{2}  \left(\frac{|\lambda-1|}{2}\right)^{\frac{1}{1+\alpha_k}}\cdot  d_A(q,   q_j).
\end{align*}
In any case, if
$|\lambda-1|\ge 1$,  then
\b{equation}\label{e14}
d_A(w_\lambda, w_{\lambda,j})\ge    \frac{1}{2}  \left(\frac{|\lambda-1|}{2}\right)^{\frac{1}{1+\alpha_k}}\cdot  d_A(q,   q_j).
\end{equation}

Let    $b>0$   be   a constant   such that
for every sequence $r_j\ra 0$,  the inequality
$d_A(w_\lambda, w_{\lambda, j})\ge b \cdot d_A(q, q_j)$    holds for all  sufficiently large $j$.  
     Let $\tilde w_{\lambda, j}= r^{1+\beta_l}_j e  +  F(w_\lambda)\in r^{1+\beta_l}_j e+ L'$.    By  the quasisymmetric condition and  (\ref{e9})
 \b{align*}
\eta_1(1)  \cdot \frac{l_{F^{-1}}( F(w_\lambda),   r_j)}{r_j} & \ge
 \frac{d_A(w_\lambda,\,  F^{-1}(\tilde w_{\lambda, j}))}{r_j} \\
 & \ge  \frac{d_A(w_\lambda, \, w_{\lambda, j})}{r_j}\\
 &\ge  b   \cdot \frac{d_A(q, q_j)}{r_j}\ge  \frac{b}{\eta_1(1)} \cdot  \frac{l_{F^{-1}}(q',  r_j)}{r_j}.
\end{align*}
  Hence
 $$\liminf_{j\ra\infty}  \frac{l_{F^{-1}}( F(w_\lambda),   r_j)     }{r_j}  \ge \frac{b}{(\eta_1(1))^2} \cdot \liminf_{j\ra\infty} \frac{l_{F^{-1}}(q',  r_j)}{r_j}\ge
 \frac{b}{(\eta_1(1))^2} \cdot  l_{F^{-1}}(q').$$
Since this holds for   every   sequence $r_j\ra 0$, we have
$l_{F^{-1}}(F(w_\lambda))\ge \frac{b}{(\eta_1(1))^2} \cdot   l_{F^{-1}}(q')$.
  Therefore,
  \b{equation}\label{e15}
 \text{L}_F(w_\lambda)\le  \frac{(\eta_1(1))^2}{b}\cdot  \text{L}_F(q).
\end{equation}

Combining  (\ref{e14}) and  (\ref{e15})   we  see that the following holds for all $|\lambda-1|\ge 1$:
$$ \text{L}_F(w_\lambda)\le  2{(\eta_1(1))^2}\left(\frac{2}{|\lambda-1|}\right)^{\frac{1}{1+\alpha_k}}\cdot  \text{L}_F(q).$$


\end{proof}

Recall the grading   ${\mathcal H}_n=V_A\oplus \mathcal Z({\mathcal H}_n)$,  where $V_A=U_1\oplus \cdots \oplus U_k$.
   Let $\pi_1:   {\mathcal H}_n\ra V_A$  be the projection with respect to the above grading.
  Let $L$ be a left coset of $U_1$.    Notice that the restriction $\pi_1|_L$ is injective  and $\pi_1(L)$ is an affine subspace of $V_A$.
  A subset $H\subset L$ is called a hyperplane of   $L$  if $\pi_1(H)$ is a hyperplane of $\pi_1(L)$.  Similarly,
  a subset $A\subset L$ is called a line in $L$ if  $\pi_1(A)$ is a  line in $\pi_1(L)$.

\b{Le}\label{halfspace}
Let $L$  be   a left coset of $U_1$.
   Suppose $p,q\in L$ are such that
 $l_F(p)> (C_2)^{2m} \text{L}_F(q)$,  where    
 $$C_2=\max\{102\eta_1(1)\eta(1),  3\eta(1)(\eta_1(1))^3\}$$
  and    $m=\dim(U_1)$.
    Then  there is a hyperplane $H$   of $L$  passing through $q$  and one
   component $H_-$ of $L\backslash H$ such that $l_F(x)\le  (C_2)^{2m} \text{L}_F(q)$  for all $x\in H_-$.

\end{Le}

\b{proof}
Let $S$ denote the space of directions of $L$ at $q$.   We shall define two subsets  $G$, $B$ of $S$.
 A point $s\in S$ lies in $G$ if  $\text{L}_F(x)\le
 C_2^m \text{L}_F(q)$  for  every $x\not=q$  in the direction of $s$.
   A point $s\in S$ lies in $B$ if  $l_F(x)>
 C_2^{2m} \text{L}_F(q)$  for    some  $x\not=q$  in the direction of $s$.
  Clearly $G\cap B=\emptyset$.    Let $s_1\in S$ be the direction of $p$,  and $s_2\in S$ the point in $S$ opposite to $s_1$.    Then $s_1\in B$ since
$l_F(p)> C_2^{2m} \text{L}_F(q)$.  Lemma \ref{key} implies $\text{L}_F(x)\le 3(\eta_1(1))^2 \text{L}_F(q)$ for any point $x\not=q$  such that $q\in xp$.  Hence $s_2\in   G$.

Let $H(B)\subset S$ be the convex hull of $B$ in the sphere $S$.    Then  for any $y\in H(B)$,
  there are $m$ points $x_1, \cdots, x_m\in B$ such that $y$ lies in the spherical   simplex  $\Delta_1$  spanned by
   $x_1, \cdots, x_m$.    Let $\Delta_i$ be the spherical simplex spanned by $x_i, \cdots, x_m$.
  Then there are $y_i\in \Delta_{i}$  with   $y_1=y$  such that  $y_i\in x_i y_{i+1}$.
   Since $x_i\in B$, there exists a point $p_i\not=q$  in the direction of $x_i$ such that
   $l_F(p_i)>
 (C_2)^{2m} \text{L}_F(q)$.
  Let $q_{m-1}$ be the unique point in the direction of $y_{m-1}$ such that $q_{m-1}\in p_{m-1}p_m$.
    Inductively, let $q_i$ be the unique point in the direction of $y_i$ such that $q_i\in p_iq_{i+1}$.

 We claim $\text{L}_F(x)> C_2^{2m-1} \text{L}_F(q)$ for every $x\in  p_{m-1}p_m$; in particular,
   $\text{L}_F(q_{m-1})> C_2^{2m-1} \text{L}_F(q)$.
  Suppose not. Then
$\text{L}_F(x)\le  C_2^{2m-1} \text{L}_F(q)$ for   some  $x\in  p_{m-1}p_m$.
  Since $l_F(p_{m-1})> C_2^{2m} \text{L}_F(q)$,  we have
 $  l_F(p_{m-1})>  C_2  \text{L}_F(x)$.  Now Lemma \ref{key} implies
   $\text{L}_F(p_m)\le 3 (\eta_1(1))^2 \text{L}_F(x) \le   C_2^{2m} \text{L}_F(q)$, contradicting  $l_F(p_m)> C_2^{2m} \text{L}_F(q)$.  By   considering $q_i\in q_{i+1}p_i$  and using Lemma \ref{key}  one inductively  proves  that
 $\text{L}_F(q_{i})> C_2^{m+i} \text{L}_F(q)$.
   In particular,   $\text{L}_F(q_{1})> C_2^{m+1} \text{L}_F(q)$.
 Since $q_1$ is in the direction   of   $y_1=y$, we see that $H(B)\cap G=\emptyset$.

Now $H(B)$ is a    non-empty convex subset of the sphere $S$ and its complement is non-empty.
  It follows that there is an   open hemisphere in its complement.  Hence there is an open hemisphere in  the complement of $B$.   Now the Lemma follows from the definition of $B$.

\end{proof}

\b{Le}\label{nonzero}
 Suppose $\dim(U_1)\ge 2$.   Then
  for any bounded subset  $X\subset L$,  there exist two positive constants $M_1, M_2$ such that
  $\text{L}_F(x)\ge M_1$ and $l_F(x)\le M_2$ for all $x\in X$.

\end{Le}

\b{proof}  Let $X$  be a bounded subset of $L$.
We first show that there is some $M_1>0$ such that
   $\text{L}_F(x)\ge M_1$  for all $x\in X$.
 Suppose there is a sequence of points  $x_i\in X$ such that
$\text{L}_F(x_i)\ra 0$.   Fix a point $p\in L$ such that $F|_L$ has  non-singular differential at $p$.
  Such a point $p$ always exists:   
  $(L, d_A)$   and $(F(L), d_B)$ are  isometric to $\R^m$; 
 since  $\dim(U_1)\ge 2$ and $F|_L$ is quasisymmetric,   $F|_L$ is a.e. differentiable and its differential is a.e. non-singular.
  The quasisymmetry condition implies  $l_F(p)>0$.
Then for all sufficiently large $i$ we have
    $l_F(p)>C^{2m}_2 \text{L}_F(x_i)$.    By Lemma \ref{halfspace},
  there is a hyperplane $H_i$   passing through $x_i$  and a component $H_{i,-}$ of $L\backslash H_i$   such that
  $l_F(x)\le C_2^{2m} \text{L}_F(x_i)$ for all $x\in H_{i,-}$.
  Since the sequence $x_i$ is bounded,   a subsequence
  $H_{i_j,-}$
of  the half spaces
$H_{i,-}$  converges to an open half space $H_-$.
  Since every $x\in H_-$ lies in  $H_{i_j,-}$  for all sufficiently large $j$  and
$\text{L}_F(x_i)\ra 0$, it follows that $l_F(x)=0$ for all $x\in H_-$.  Since
$F|_L: L\ra F(L)$
is a.e. differentiable,  we  see that $F|_L$
 has    zero differential   a.e. on the open set  $H_-$  of $L$,   which is impossible.

   As a quasisymmetric map,    $F|_L: L\ra F(L)$  maps bounded sets to bounded sets. So
 $F(X)$ is bounded.   Now the first claim applied to $F^{-1}$    yields that
   there is a positive lower bound for $\text{L}_{F^{-1}}$  on $F(X)$.
  Now   (\ref{e5}) implies  that  there is a positive upper bound  for $l_F$  on $X$.

\end{proof}

It is clear from the definition of $d_A$ and $d_B$ that  lines in the left cosets of $U_1$ and $W_1$ are rectifiable  (recall we first normalized so that $\alpha_1=\beta_1=1$).

\b{Le}\label{quasisimi}
For each left coset $L$  of $U_1$, there is some constant $C>0$ such that
   $F|_L $ is a $(C^{2m+2}_2, C)$-quasisimilarity,
  where   $m=\dim(U_1)$ and $C_2$ is the constant in Lemma \ref{halfspace}.


\end{Le}

\b{proof}
First   consider the case when $m=1$.
  Lemma \ref{rec}  and the comment before Lemma \ref{quasisimi}   imply that the left cosets of $U_1$ are the only rectifiable curves in $(H_n, d_A)$. Similarly the
  left cosets of $W_1$ are    the only rectifiable curves in $(H_n, d_B)$.
 By the main result of \cite{BKR},   $F$ is absolutely continuous on a.e. left coset of $U_1$.   Let $L$ be such a left coset.
 Since   $F|_L: (L, d_A)\ra (F(L), d_B) $  is a  homeomorphism between lines (with the Euclidean metric),  it is differentiable a.e.   As $F|_L$ is absolutely continuous,
it suffices   to  bound the differential  in order to show that $F|_L $ is a quasisimilarity.
  We shall show that  $l_F(p)\le  C_1\cdot \text{L}_F(q)$ for any $p,q\in L$, where $C_1=102\eta_1(1)$.
Suppose there are two points $p, q\in L$ such that
  $l_F(p)> C_1\cdot \text{L}_F(q)$.     By Lemma \ref{key},
   $\text{L}_F(x)\ra 0$ as $d_A(p, x)\ra \infty$ ($x\in L$).  This implies
  $\text{L}_{F^{-1}}(y)\ge l_{F^{-1}}(y) \ra \infty$ as $d_B(y, F(p))\ra \infty$ ($y\in F(L)$).
  However,    $l_F(p)> C_1 \text{L}_F(q)$  implies
 $l_{F^{-1}}(F(q))> C_1 \cdot \text{L}_{F^{-1}}(F(p))$.  By Lemma \ref{key} again we obtain
  $\text{L}_{F^{-1}}(y)\ra  0$, which is a contradiction.

From now on we assume    $m\ge  2$.  Denote $L'=F(L)$.  In this case, both $F|_L$ and $F^{-1}|_{L'}$ have the following properties:
   (1)  absolutely continuous,
 (2)  differentiable almost everywhere and the differential is almost everywhere nonsingular,
   (3)    absolutely continuous on almost all curves.
  It follows that to show $F|_L$   is a  $(C^{2m+2}_2, C)$ quasisimilarity, it suffices to show that
    there is    a  set of full measure $E\subset L$ such that
    $l_F(x)\le 3 (\eta_1(1))^2 C^{2m+1}_2 \text{L}_F(y)$  for all $x, y\in  E$.
    We shall prove   by contradiction.  So suppose the above statement is not true.
 Then   in particular  there are two points  $p, q\in L$ such that
   $l_F(p)> 3  (\eta_1(1))^2  C^{2m+1}_2 \text{L}_F(q)$.

 We observe that it suffices to show that there is a  constant $b_0>0$ such that
$l_F(x)\le b_0$  for all $x$ in a full measure subset of $L$:
  the condition   $l_F(p)> 3  (\eta_1(1))^2  C^{2m+1}_2 \text{L}_F(q)$
  implies that   $l_{F^{-1}}(F(q))> 3  (\eta_1(1))^2  C^{2m+1}_2 \text{L}_{F^{-1}}(F(p))$.
  Then Lemma \ref{key} implies
  $\text{L}_{F^{-1}}(y)\ra 0$ as $y\in L'$ goes to infinity along the line through $F(p)$  and $F(q)$.
  Fix a point $y_0$ such that
  $$\text{L}_{F^{-1}}(y_0)<\min\left\{\frac{1}{b_0  \eta_1(1) C_2^{2m}}, \; \frac{l_{F^{-1}}(F(q))}{C_2^{2m}}\right\}.$$
 By Lemma \ref{halfspace},
  there is a hyperplane   $H'$ passing trough $y_0$  and a component $H'_-$  of $L'\backslash H'$ such that
     $ l_{F^{-1}}(y)<  \frac{1}{b_0\eta_1(1)}$  for all $y\in H'_-$.
  Since $F^{-1}$ is differentiable a.e. on $L'$ and is $\eta_1$-quasisymmetric,
  we have $ \text{L}_{F^{-1}}(y)<  \frac{1}{b_0}$  for   a.e. $y\in H'_-$.
It follows that
   $l_F(x)> b_0$ for a.e.  $x\in F^{-1}(H'_-)$, contradicting  the assumption that
$l_F(x)\le b_0$  for all $x$ in a full measure subset of $L$.

We next show that $l_F(x)$ is essentially bounded on $L$.

Since
$l_F(p)> 3  (\eta_1(1))^2  C^{2m+1}_2 \text{L}_F(q)$,
  by Lemma   \ref{halfspace},
    there is a hyperplane $H_1$ passing through $q$  and one component $H_{1,-}$ of $L\backslash H_1$ such that 
$$l_F(x)\le  (C_2)^{2m} \text{L}_F(q)<\frac{1}{3(\eta_1(1))^2 C_2}  l_F(p)$$
 for all $x\in H_{1,-}.$  The quasisymmetry condition then implies
 $$L_F(x)\le \eta(1) l_F(x)\le \frac{\eta(1)}{3(\eta_1(1))^2 C_2}  l_F(p)$$
 for   a.e.  $x\in H_{1,-}.$ 
    Let $\tau:  L\ra L$ be the geodesic symmetry about $p$,  that is,  for any $x\in L$,  $\tau(x)$ is such that $p$ is the
    midpoint of $x\tau(x)$.    Now  Lemma \ref{key}  implies that   for a.e. $y\in \tau(H_{1,-})$  we have
    $$\text{L}_F(y)\le 3 (\eta_1(1))^2 \text{L}_F(\tau(y))\le  3 \eta(1)(\eta_1(1))^2  C_2^{2m} \text{L}_F(q). $$
  If $p\in H_1$, then we are done since   now $\text{L}_F(x)$ (and hence $l_F(x)$)  is bounded on a  full measure subset  of 
 $L\backslash H_1=H_1\cup  \tau(H_1)$.    So we assume $p\notin H_1$.
Let $B_1$ be the part of a cylinder in $L$  between $H_1$ and $\tau(H_1)$  with center line  passing through  $p$  and
    perpendicular to the hyperplane $H_1$.     Then $B_1$ is bounded.  By Lemma \ref{nonzero}
      there are positive numbers $M_1, M_2$  such that
   $\text{L}_F(x)\ge M_1$ and $l_F(x)\le M_2$ for all $x\in B_1$.

  Since  we assume
    $l_F(x)$ is not essentially bounded,   there is some hyperplane $\tilde{H}_1$  parallel to $H_1$
 such that (1) $F|_L$ is differentiable at some $q_1\in B_1\cap \tilde{H}_1$;  (2)  there is some $p_1\in \tilde{H}_1$   with  $l_F(p_1)>C^{2m}_2\cdot \eta(1) M_2$.  Since $F|_L$ is differentiable at
 $q_1$,  we have $\text{L}_F(q_1)\le \eta(1) l_F(q_1)\le \eta(1) M_2$.    So
  $l_F(p_1)>C^{2m}_2\cdot   \text{L}_F(q_1)$.  Now Lemma \ref{halfspace}  and Lemma \ref{key}  imply
    that there
  is a hyperplane $H_2$ passing through $q_1$ and a component $H_{2, -}$  of $L\backslash H_2$ such that
   $\text{L}_F(x)$  is  essentially  bounded   from above  on
$H_{2, -}$    and $\tau_1(H_{2, -})$, where $\tau_1$ is the geodesic symmetry about $p_1$.
  If $p_1\in H_2$, then we are done as indicated above. So we assume $p_1\notin H_2$. In this case,  $H_1$ and $H_2$ are not parallel.    We proceed inductively and eventually find $m$ hyperplanes $H_1$,
 $H_2$,  $\cdots$,  $H_m$,  $m$  half spaces $H_{i, -}$  and points $p_0=p, p_1, \cdots,  p_{m-1}$ with the following properties:
  (1)   $\text{L}_F(x)$  and hence  $l_F(x)$ is  essentially  bounded   from above  on the union
  $Q:=\cup_i H_{i, -}\cup \cup_i \tau_{i-1}(H_{i, -})$, where $\tau_{i-1}$  ($\tau_0=\tau$) is the geodesic symmetry about the point $p_{i-1}$;
  (2)   The complement of $Q$ in $L$ is compact.
  By Lemma \ref{nonzero},  $l_F(x)$ is uniformly bounded on $L\backslash Q$. It follows that
  $l_F(x)$  is essentially   bounded on $L$,  and we are done.

\end{proof}

\section{Proof of the main Theorems}\label{proofs}

In this Section we finish the proofs of the theorems in the Introduction.

We use the notation from Section \ref{leaf},  see the paragraphs before Lemma \ref{key}.

Notice that
$\bigoplus_{i=1}^{k-1}U_i\oplus \mathcal Z({\mathcal H}_n)$  is a  Lie subalgebra (actually an ideal) of   ${\mathcal H}_n$.
 Let  $H$ denote   
    the corresponding connected Lie subgroup of $H_n$.

\b{Le}\label{para}
Suppose $k\ge 2$.     Then  two left cosets of $U_1$  lie in the same left coset of   $H$  
 if and only if the Hausdorff distance between them  is finite.

\end{Le}

\b{proof}
Let $L_1$, $L_2$ be two left cosets of $U_1$.
   After  applying   a  left translation, we may assume
  $L_1=U_1$ and $L_2=g*U_1$ for some  $g=x_1+\cdots +x_k+x_{k+1}\in \mathcal H_n$  with  $x_i\in U_i$.
  Note that $L_1$ and $L_2$ lie in the same left coset of $H$ if and only if $x_k=0$.

First assume $L_1$ and $L_2$ lie in the same left coset of $H$.  Then $x_k=0$.
   For $t\in U_1$,  Lemma \ref{structure} (2) implies
 $$(-t)*g* t=(-t)*(x_1+\cdots  + x_{k-1} +x_{k+1})*t=x_1+\cdots  + x_{k-1}+x_{k+1}.   $$
   We see that $d_A(t, g*t)=||(-t)*g*t||_A$
 is independent of $t\in U_1$.
  Hence  the Hausdorff distance between   $L_1$ and $L_2$    is finite.

Next we assume $L_1$ and $L_2$ lie in    distinct   left cosets of $H$.  Then $x_k\not=0$.
  There exists $v\in U_1$ such that $[x_k, v]\not=0$.    Let $t_1=a v$ with $a\in \R$.
 Let $t_2\in U_1$.  We consider $d_A(t_1, \;  g*t_2)$.  Calculate
  $$(-av)*g*t_2=\left((x_1+t_2-av)+x_2+\cdots +x_k+(x_{k+1}+\frac{1}{2}[x_k,  av+t_2])\right). $$
Suppose the Hausdorff distance between $L_1$ and $L_2$ is finite. Then there is  some constant $C>0$ such that,
   for any $a\in \R$,    there is some $t_2\in U_1$  satisfying
  \b{align*}
 C\ge d_A(av, g*t_2) & =||(-av)*g*t_2||_A  \\
 &=|x_1+t_2-av|+\sum_{i=2}^k |x_i|^{\frac{1}{\alpha_i}}+\left|x_{k+1}+\frac{1}{2}[x_k,  av+t_2]\right|^{\frac{1}{1+\alpha_k}}.
\end{align*}
   Let $u=t_2-av$.  Then $u$ is uniformly bounded and $t_2=av+u$.
    It follows that $[x_k, u]$ is uniformly bounded.  As $[x_k, v]$ is fixed and nonzero,  we see that
$[x_k, av+t_2]=[x_k, u]+2a[x_k, v]$  is unbounded as $a\ra \infty$.   The contradiction shows that the Hausdorff distance between
 $L_1$ and $L_2$ is infinite.

\end{proof}

\b{Le}\label{housk}
For any $x_k, x'_k\in U_k$ and any $h\in H$, we have
 $$d_A(x_k*H, x'_k*H)=d_A(x_k*h,  x'_k*H)=|x'_k-x_k|^{\frac{1}{\alpha_k}}.$$

\end{Le}

\b{proof}
 Write $h=x_1+\cdots +x_{k-1}+x_{k+1}$.    For any $h'=x'_1+\cdots+ x'_{k-1}+x'_{k+1}\in H$,
 we have
\b{align*}
  d_A(x_k*h,   x'_k*h') &=||(-h)*(-x_k)*x'_k*h'||_A\\
&=\left|\left| (x'_1-x_1)+\cdots  +(x'_{k-1}-x_{k-1})+(x'_k-x_k)+  (x'_{k+1}-x_{k+1}+E) \right| \right|_A,
\end{align*}
    where $$E=\frac{1}{2}[x'_k-x_k,  x_1+x'_1]+\frac{1}{2}
   \left[-\sum_{i=1}^{k-1} x_i,  \sum_{i=1}^{k-1} x'_i\right].$$
  Now it is clear that $d_A(x_k*h,   x'_k*h')\ge |x'_k-x_k|^{\frac{1}{\alpha_k}}$ for any $h'\in  H$.
Furthermore,   the equality holds
   for $h'=x_1+\cdots +x_{k-1}+(x_{k+1}-[x'_k-x_k, x_1])$.  
  The Lemma follows.

\end{proof}

The set of left cosets of $H$ in $H_n$ can be identified with $U_k$ via $x_k\ra x_k*H$. Lemma \ref{housk}   implies that this  set equipped with
   the minimal distance is isometric to $(U_k, |\cdot|^{\frac{1}{\alpha_k}})$.

  Set
$K=\bigoplus_{i=2}^{k-1}U_i\oplus \mathcal Z({\mathcal H}_n)$.
A similar (and easier) calculation as in the proof of Lemma \ref{housk}   yields the following:

\b{Le}\label{hous1}
    We have
 $d_A(g*U_1, \,  g'*U_1)=d_A(g*x, \, g'*U_1)=d_A(g, g')$
     for any $g, g'\in  K$ and any $x\in U_1$.

\end{Le}

 Lemma \ref{hous1}  implies that the set of left cosets of $U_1$ in $H$ can be identified with $K$,  and this set equipped
 with the minimal distance is isometric to $(K, d_A)$.

The proof of the following Lemma  is almost the same as that of   Lemma 3.9 in \cite{X3}.
 So we omit the proof here. The main point  is  that different left cosets diverge sublinearly.

\b{Le}\label{bilip1}
  Suppose $\alpha_1=\beta_1=1$.
 Then there is a constant $C$ such that for any left coset $L$  of $U_1$,  $F|_L$  is
 a   $(K, C)$ quasisimilarity, where $K$ depends only on $\eta$.

\end{Le}


\noindent
{\bf{Proof of Theorem  \ref{main3}}.}
Let $A, B$   be diagonalizable    derivations  with positive eigenvalues.
       Suppose   $k\ge 2$.
Let $F:   (H_n, d_A)  \ra (H_n, d_B) $  be a  quasisymmetric map.
   We shall use the notation from Section \ref{leaf}.
Then     $F:   (H_n, d_1)  \ra (H_n, d_2) $
 is  $\eta$-quasisymmetric
  for some $\eta$,  where $d_1=d_A^{{\alpha_1}}$  and $d_2=d_B^{{\beta_1}}$  .
  Lemma \ref{bilip1}  implies that
 there is a constant $C$ such that for any left coset $L$  of $U_1$,  $F|_L:  (L, d_1)  \ra  (F(L),  d_2) $  is  a
  $(K, C)$ quasisimilarity, where $K$ depends only on $\eta$.

   Let $p, q\in H_n$ be arbitrary. If they lie in the same left coset $L$ of $U_1$, then
   $d_2(F(p), F(q))\le CK d_1(p,q)$.
  Now suppose $p\in L_1$, $q\in L_2$.  Pick $x\in L_1$ such that $d_1(p, x)=d_1(p, q)$.
   Then
  $$d_2(F(p), F(q))\le \eta(1)\cdot  d_2(F(p), F(x))\le \eta(1) CK d_1(p, x)=\eta(1)CK d_1(p,q).$$
  So we have an upper bound for $d_2(F(p), F(q))$.    The same  argument applied to $F^{-1}$   yields
    a lower
  bound for $d_2(F(p), F(q))$.   Hence $F$ is  biLipschitz.   The theorem then follows.

\qed

\noindent
{\bf{Proof of Theorem  \ref{main4}}.}
First we suppose $A$ and $B$ have the same invariants, that is,
   $l=k$, $\text{dim}(U_i)=\text{dim}(W_i)$   and there is    some $\lambda>0$ such that $\beta_i=\lambda \alpha_i$.   We need to show that $(H_n, d_A)$ and $(H_n, d_B)$ are quasisymmetric.
   After replacing $d_A$ with $d_A^{\alpha_1}$  and
$d_B$ with $d_B^{\alpha_1}$,
    we may assume  $\alpha_1=\beta_1=1$. Then $\beta_i=\alpha_i$.   Fix some $e\in   \mathcal Z({\mathcal H}_n)\backslash\{0\}$.
  By Lemma   \ref{structure},  for $1\le i<(k+1)/2$,
there is a basis
  $e_1, \cdots, e_{m_i}$  for $U_i$ and a basis $\eta_1, \cdots, \eta_{m_i}$ for $U_{k+1-i}$ such that
   $[e_s, \eta_t]=\delta_{st} e$;   if  $i=(k+1)/2$,   then $m_i=2k_i$ is even   and  there is a basis
 $e_1,  \eta_1, \cdots,  e_{k_i}, \eta_{k_i}$ of  $U_i$  such that  $[e_s, \eta_t]=\delta_{st} e$, $[e_s, e_t]=[\eta_s, \eta_t]=0$.
  Similarly,
  for $1\le i< (k+1)/2$,
there is a basis $e'_1, \cdots, e'_{m_i}$  for $W_i$ and a basis $\eta'_1, \cdots, \eta'_{m_i}$ for $W_{k+1-i}$ such that
   $[e'_s, \eta'_t]=\delta_{st} e$;  and if   $i=(k+1)/2$,   then $m_i=2k_i$ is even   and  there is a basis
 $e'_1,  \eta'_1, \cdots,  e'_{k_i}, \eta'_{k_i}$  for $W_i$
 such that  $[e'_s, \eta'_t]=\delta_{st} e$, $[e'_s, e'_t]=[\eta'_s, \eta'_t]=0$.
 Now define a map $G: V_A \ra V_B$ as follows.  For each $i<(k+1)/2$, $G|_{U_i}$ is given by
$G(\sum_s x_s e_s)=\sum_s x_s e'_s$;
  if $i>(k+1)/2$,   define  $G(\sum_s x_s \eta_s)=\sum_s x_s \eta'_s$; and if
 $i=(k+1)/2$, define 
$$G(\sum_s (x_s e_s+y_s \eta_s))=\sum_s (x_s e'_s+y_s \eta'_s).$$
  Define $F:(H_n, d_A)\ra (H_n, d_B)$ by   $F(x+x_{k+1})=G(x)+x_{k+1}$, where $x\in V_A$ and $x_{k+1}\in \mathcal Z({\mathcal H}_n)$.    It is now easy to check that $F$ is an isometry.

Conversely,
let  $F: (H_n,  d_A)\ra (H_n,   d_B)$  be a quasisymmetry.
   By Proposition \ref{prefo},      $k\ge 2$ if and only if $l\ge 2$; furthermore, in this case,  $\dim(U_1)=\dim(W_1)$  and 
  $F$   maps left cosets of $U_1$ to left cosets of $W_1$.
  The  conclusion of the theorem  clearly    holds if $k=l=1$. So from now on we shall assume $k\ge 2$ and $l\ge 2$.
 After replacing $d_A$ with $d_A^{\alpha_1}$  and
$d_B$ with $d_B^{\alpha_1}$,
   we may assume $\alpha_1=\beta_1=1$.
By   Theorem \ref{main3},  $F$ is   biLipschitz.
Lemma \ref{para}  says two left cosets $L_1$ and $L_2$  of $U_1$  lie in the same left coset of
 $H:=\bigoplus_{i=1}^{k-1} U_i \oplus \mathcal Z({\mathcal H}_n)$ if and only if the   Hausdorff distance between them is finite.   The same is true for left cosets of $W_1$
  and $\tilde H: =\bigoplus_{i=1}^{l-1} W_i \oplus \mathcal Z({\mathcal H}_n)$.   It follows that $F$ maps left cosets of $H$ to left cosets of $\tilde H$.
Now Lemma \ref{housk}  and the remark after that Lemma imply that $F$ induces a biLipschitz map from
 $(U_k,  |\cdot|^{\frac{1}{\alpha_k}}) $ to  $(W_l,  |\cdot|^{\frac{1}{\beta_l}}) $.  From this we conclude that
 $\dim(U_k)=\dim(W_l)$ and $\alpha_k=\beta_l$.

 Now we consider the restriction of $F$ to a left coset of $H$, which is biLipschitz. This can be viewed as a biLipschitz map from
 $(H, d_A) $ to $(\tilde H, d_B)$.  We also know that $F$ maps left cosets of $U_1$ to left cosets of $W_1$.
Now Lemma \ref{hous1}  and the remark after that Lemma imply that $F$ induces a biLipschitz map from
  $(K, d_A)$ to   $(\tilde K, d_B)$, where $K=\bigoplus_{i=2}^{k-1} U_i \oplus \mathcal Z({\mathcal H}_n)$  and
$\tilde K=\bigoplus_{i=2}^{l-1} W_i \oplus \mathcal Z({\mathcal H}_n)$.    Now an induction argument finishes the proof of Theorem \ref{main4}.

\qed

Theorem \ref{main2} follows from Theorem \ref{main4}   since  $G_A$ and $G_B$ are quasiisometric if and only if
 $(H_n, d_A)$ and $(H_n, d_B)$ are quasisymmetric.
  Theorem \ref{main1} follows from Theorem \ref{main3}  since any quasiisometry $f: G_A\ra G_B$ induces a boundary map
 $\p f: \p G_A\ra \p G_B$, which is a quasisymmetric map,  and  $f$ is a   almost  similarity if and only if $\p f$ is biLipschitz (after possibly snowflaking the metric
 $d_B$).
 For more details on these implications, the reader is referred to \cite{SX}.

 \addcontentsline{toc}{subsection}{References}

\noindent Address:

\noindent Xiangdong Xie: Dept. of Mathematics  and   Statistics,   Bowling Green  State  University,
  Bowling Green,  OH,   U.S.A.\hskip .4cm E-mail:   xiex@bgsu.edu

\end{document}